\documentclass{article}
\usepackage{amsmath}
\usepackage{amssymb}
\usepackage{graphicx}
\usepackage{hyperref}
\usepackage{tikz}

\newlength\Colsep
\title{Parallel Magnetic Resonance Imaging}

\author{Martin Uecker}

\begin{document}
\maketitle

\section{Introduction}

The main disadvantages of Magnetic Resonance Imaging (MRI)
are its long scan times and, in consequence, its sensitivity to motion.
Exploiting the complementary information from 
multiple receive coils, parallel imaging is able to recover images 
from under-sampled k-space data and to accelerate the
measurement~\cite{Hutchinson88,Kelton89,Kwiat91,Carlson93,Ra93,Sodickson97,Pruessmann99}. 
Because parallel magnetic resonance imaging can be
used to accelerate basically any imaging sequence it
has many important applications. Parallel imaging
brought a fundamental shift in image reconstruction:
Image reconstruction changed from a simple direct Fourier 
transform to the solution of 
an ill-conditioned inverse problem. 
This work gives an overview of image reconstruction
from the perspective of inverse problems. After 
introducing basic concepts such as regularization, discretization, 
and iterative  reconstruction, advanced topics are 
discussed including algorithms for auto-calibration,
the connection to approximation theory, and the combination 
with compressed sensing.

\section{Parallel Imaging as Inverse Problem}

\subsection{Forward Model}

\begin{figure}
\begin{center}
\includegraphics[width=\textwidth]{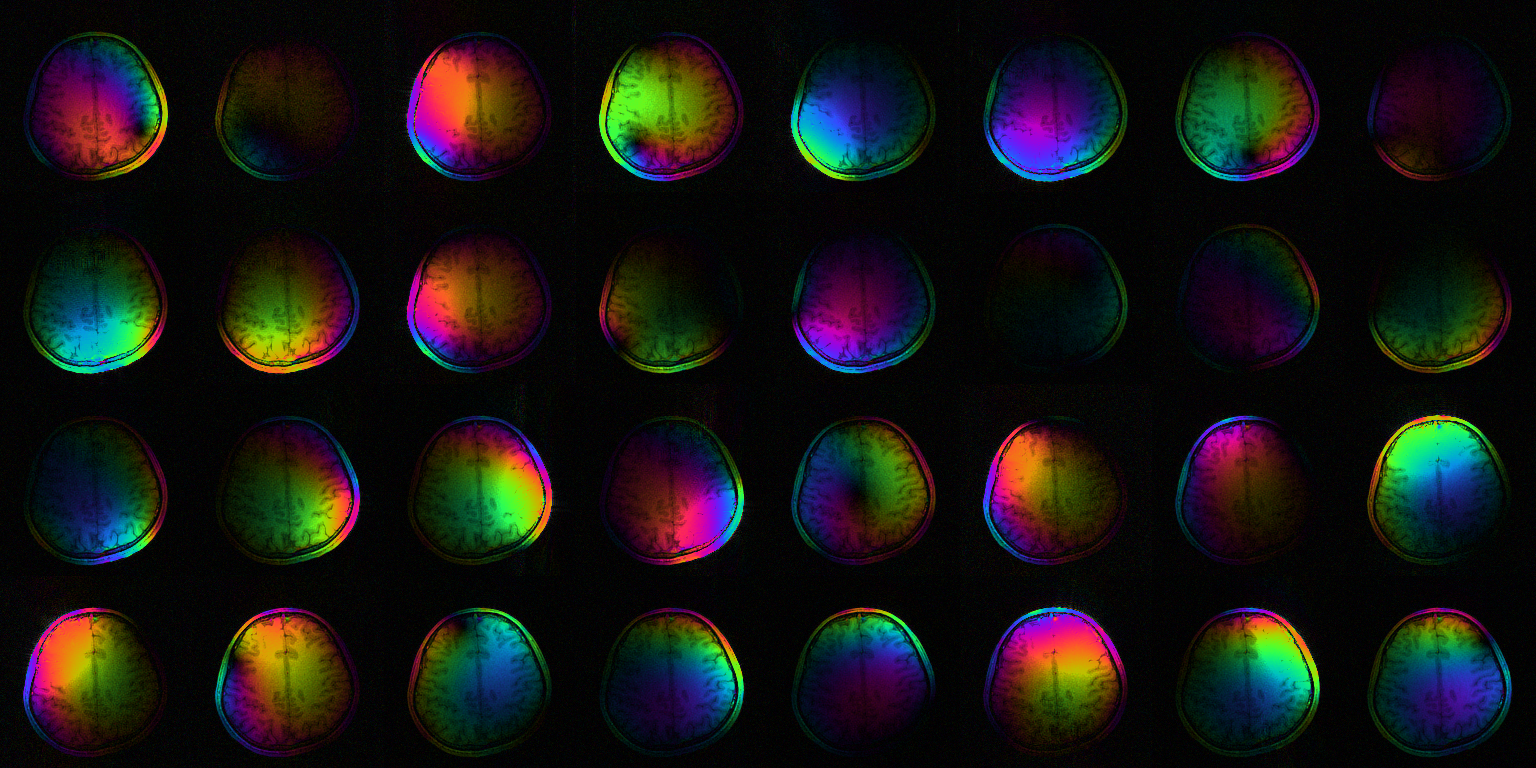}
\end{center}
\caption{All individual coil images $m c_j$ for a phased-array coil with 32 elements (channels).
For each coil element $j = 1, \dots, 32$, the complex-valued magnetization image $m$ 
is modulated by its unique receive sensitivity $c_j$. The phase is color coded.}\label{fig:cimgs}
\end{figure}

\begin{figure}
\begin{center}
\includegraphics[width=\textwidth]{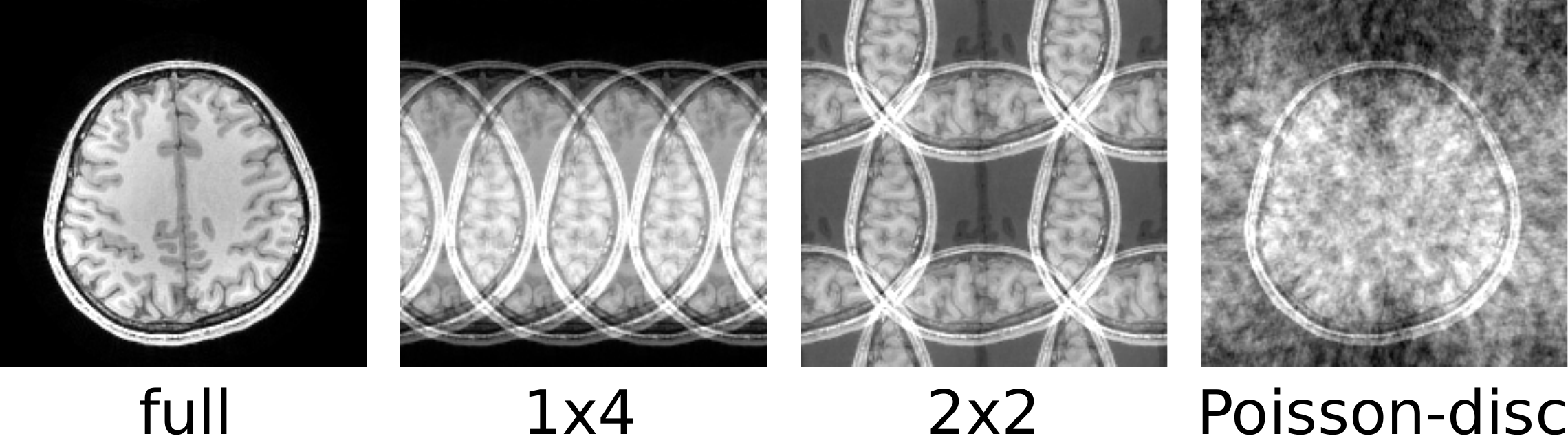}
\end{center}
\caption{Under-sampling in k-space causes aliasing artifacts in the image domain.
Images of a human brain were reconstructed with a direct Fourier transform
from fully-sampled data (full) and four-fold
under-sampled data (regular under-sampling in one dimension (1x4), two 
dimensions (2x2), and Poisson-disc sampling).}
\end{figure}


The signal from a receive-coil array with $N$ channels is given by~\cite{Roemer90}:
\begin{align}
	y_j(t) = \int_{\Omega} \textrm{d}\vec r \, m(\vec r) c_j(\vec r) e^{- 2\pi i \vec k(t) \cdot \vec r} \qquad 1 \leq j \leq N\label{eq:SIG}
\end{align}
The complex-valued magnetization image $m$ represents 
the state of the transverse magnetization 
of the excited spins in the field-of-view (FOV) $\Omega \subset \mathbb{R}^d$ at 
the time of image acquisition. In MRI, typically a volume ($d = 3$) or a 
thin slice ($d = 2$) of proton spins is excited using a resonant radio-frequency pulse.
The image (or volume) is modulated by the complex-valued sensitivities $c_j$
of all receive coils (Fig.~\ref{fig:cimgs}). The k-space signals $y_j(t)$ are then given by the 
Fourier transform of the coil images $c_j m$ sampled at discrete
time points $t_l$ along a given k-space trajectory $\vec k(t)$.
Equation~\ref{eq:SIG} neglects relaxation 
and off-resonance effects during the acquisition, which is possible if all 
samples at time points $t_l$ are acquired in a small window around 
the echo time (TE) after excitation. Relaxation and off-resonance effects from the excitation
of the spins until the echo time are incorporated into the image and 
define the image contrast. Because in most cases it is not possible to 
acquire all data in this short acquisition window, samples have to be
acquired in a repeated series of identical experiments, which
always restore the magnetization image to exactly the same
state.\footnote{Single-shot Echo-Planar Imaging (EPI)
and spiral imaging sequences acquire all data in a single
acquisition. These sequences are fast, but image quality is
compromised by blurring, distortions, and phase cancellation,
due to relaxation and off-resonance effects. Parallel imaging 
can be used to shorten the acquisition window and reduce these artifacts.}
This requirement to repeat the basic experiment many times
is the reason for the long scans times in MRI. The goal of parallel 
imaging is to reduce the amount of data required to reconstruct images
by optimally exploiting the complementary information from 
multiple receive coils. Although there are fundamental limits
to the encoding power of the receive sensitivities, it has
the potential to accelerate MRI by a factor of about four
in each spatial dimension~\cite{Wiesinger04}.

\subsection{Image Reconstruction}

\begin{figure}
\begin{center}
\includegraphics[width=0.7\textwidth]{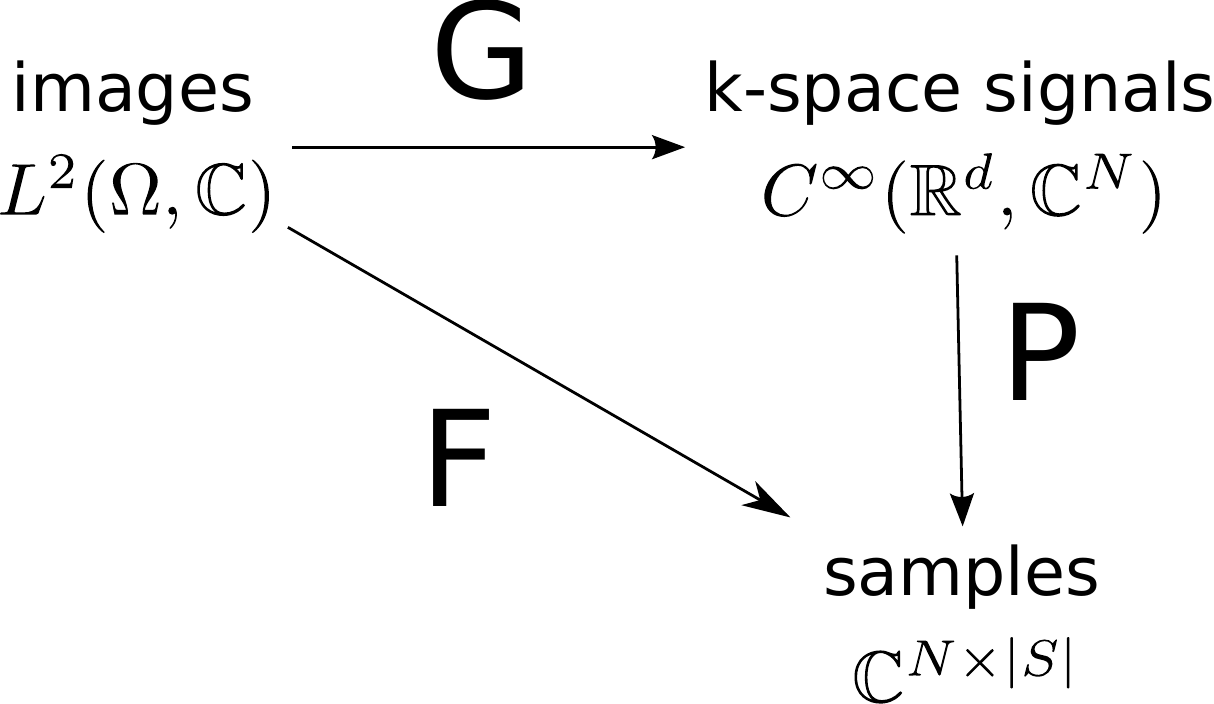}
\end{center}
\caption{Parallel imaging is an inverse problem. The forward model 
is given by a composition of a physical model and a sampling operator
$F = P \circ G$. The operator $G$ maps
magnetization images $m \in L(\Omega, \mathbb{C})$ on a field-of-view (FOV) $\Omega$
to the ideal k-space signals $f \in L(\mathbb{C}^d, \mathbb{C}^N)$ from $N$ channels.
The sampling operator $P$ maps ideal signals $f$ to k-space samples $f_j(\vec k(t_l))$ 
for all channels $j = 1, \dots, N$ and at all sample
locations $\vec k(t_l) \in S$.}\label{fig:spaces}
\end{figure}

If the coil sensitivities are known (e.g. from a pre-scan), image
reconstruction for parallel imaging can be formulated as a linear inverse 
problem with discrete data~\cite{Bertero85}.
Mathematically, the forward problem is given by
an operator $F$ which maps the magnetization image $m \in L(\Omega, \mathbb{C})$
of excited spins in a FOV $\Omega \subset R^d$ to the sample values $y_{j}(t_l)$. 
This operator can be thought of as the composition $F = P \circ G$ of a physical model $G$ 
for hypothetical multi-dimensional continuous k-space signals and a 
sampling operator $P$ (see Fig.~\ref{fig:spaces}). 
The operator $G$ is given by
\begin{align}
	G: L(\Omega, \mathbb{C}) & \rightarrow C^{\infty}(\mathbb{R}^N, \mathbb{C}) \nonumber\\
		m & \mapsto {f} \qquad \textrm{with}\quad f_j(\vec k) = \langle m, enc^{j,\vec k} \rangle~,
\end{align}
with the encoding functions defined as
$enc^{j,\vec k}(\vec r) := \overline{c_j(\vec r)} e^{+ 2\pi i \vec k(t) \cdot \vec r}$
and a scalar product defined on $L(\Omega, \mathbb{C})$ (anti-linear in the second argument).
The sampling operator $P$ evaluates the ideal k-space signals ${f_j}$
at the sample locations $\vec k(t_l) \in S$ in k-space. It is assumed 
that the sampling process corresponds to the point-evaluation of ideal k-space 
functions, i.e. $y_{j}(t_l) = f_j(\vec k(t_l))$. The sample values $f_j(S)$ are 
corrupted by additive white Gaussian noise. Although the noise
is typically correlated between receive channels, this correlation can be 
removed with a whitening step.

\begin{figure}
\begin{center}
\includegraphics[width=\textwidth]{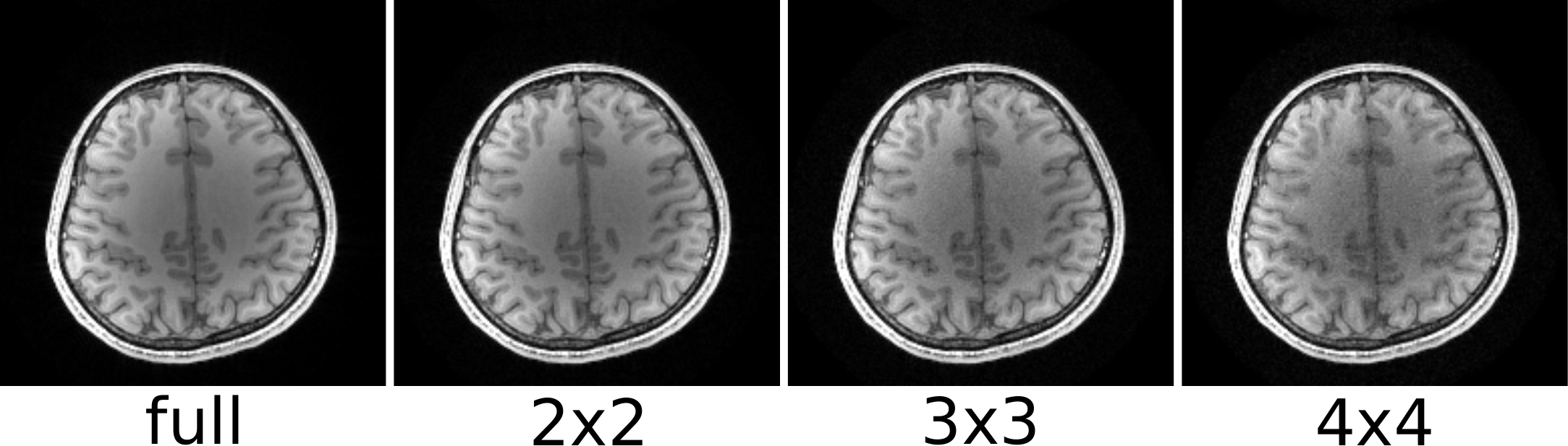}
\end{center}
\caption{Images of a human brain reconstructed using parallel imaging
from fully-sampled data (full) and under-sampled data (acceleration by 2x2, 3x3, and 4x4)
acquired with a 32-channel coil.  $l_1$-wavelet regularization was used
to suppress noise in the reconstruction.}\label{fig:reco}
\end{figure}

A variational solution to the inverse problem can be defined as the 
minimizer of a functional, i.e.
\begin{align}
	\hat m_{\alpha} := \operatorname{argmin}_x \|F x - y \|_2^2 + \alpha R(x)~.\label{eq:var}
\end{align}
The functional is composed of a least-squares data fidelity 
term (which alternatively may also include weighting or use a robust norm~\cite{Johnson12,Cheng14})
and an additional regularization term $R$.
Discretized versions of this minimization problem are the basis of SMASH
and SENSE parallel imaging methods~\cite{Ra93,Sodickson97,Pruessmann99,Sodickson01,Pruessmann01}.
For parallel MRI, this formulation has two advantages: First, arbitrary Cartesian 
or non-Cartesian sampling schemes can be used~\cite{Pruessmann01}.
Second, the regularization term can be used to introduce prior knowledge 
about the solution.\footnote{A third advantage - which conceptually goes 
beyond parallel imaging - is the possibility to extend the 
forward model to include further physical effects in model-based
reconstruction, e.g. field maps~\cite{Sutton03,Sutton04}, motion-induced
phase maps~\cite{Liu05,Uecker09a}, motion~\cite{Odille08,Cheng12},
relaxation maps~\cite{Olafsson08,Graff08,Block09},
or diffusion models~\cite{Welsh13}, etc.}
Figure~\ref{fig:reco} shows images of a human brain 
recovered from under-sampled data by numerical 
optimization of Equation~\ref{eq:var}.

\subsection{Regularization}

\begin{figure}
\begin{center}
\includegraphics[width=\textwidth]{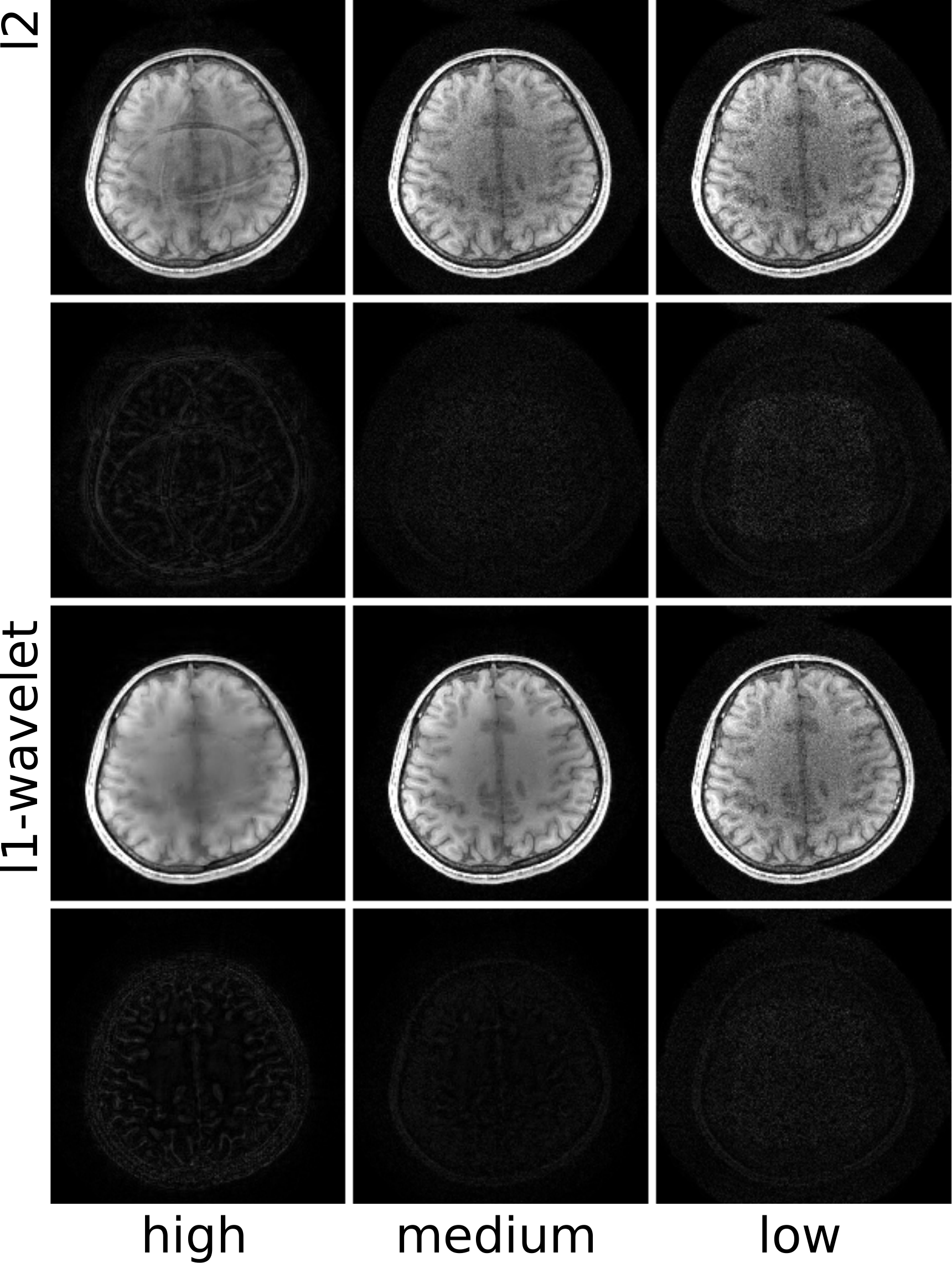}
\end{center}
\caption{Reconstruction from 16-fold under-sampled data 
using low, medium, and high amount of $l_2$ and $l_1$-wavelet 
regularization. Reconstructed images and difference images
are shown. There is a trade-off between bias (residual artifacts, blurring)
and noise depending on the amount of regularization.}\label{fig:reg}
\end{figure}

Ill-conditioning causes noise amplification during image reconstruction, which
initially limited the application of parallel imaging to only moderate
acceleration. This limitation can be overcome by incorporating prior knowledge
about the image using regularization methods~\cite{King01,Sodickson01,Lin04}.
In the simplest case, regularization may consist of
a basic quadratic penalty in the framework of a linear reconstruction,
or make use of much more sophisticated techniques which exploit
the structure of images but demand a non-linear reconstruction. For
a least-squares problem with 
quadratic regularization, i.e. $R(x) = \|\sqrt{W} (x - x_0) \|_2^2$ with a
positive definite operator $W$, the solution of Eq.~\ref{eq:var} 
is explicitly given by the formula
\begin{align}
	\hat m_{\alpha} & = x_0 + \underbrace{(F^H F + \alpha W)^{-1} F^H}_{F^{\dagger}_{\alpha}} (y - F x_0)~.\label{eq:pinv}
\end{align}
In the limit $\alpha \rightarrow 0$ this solution is called the best approximate 
and is given by the Moore-Penrose pseudo-inverse $F^{\dagger}$. It has a 
statistical interpretation  - assuming white Gaussian noise - as 
the best unbiased estimate for the image. Regularization can be interpreted
as prior knowledge and the optimizer as a maximum a posteriori (MAP) 
estimate of the image. 
Although regularization leads to a fundamental trade-off between bias and 
noise - which has to be chosen carefully for optimal image quality - it 
makes the use of higher acceleration possible. 
An optimal estimate in terms of mean squared error can only be obtained with 
regularization.

For optimal results, the prior knowledge should include as much
specific knowledge about the image as possible. For example, 
regularization can exploit smoothness in the time domain~\cite{Adluru07}, 
or exploit that changes relative to
a fixed reference image $x_0$ can be assumed to be small. The later is used 
successfully for real-time MRI~\cite{Uecker12} or dynamic contrast 
enhanced (DCE) MRI~\cite{Xu13}.
While $l_2$-regularization is simple to implement and already a clear
improvement compared to unregularized parallel imaging, much
better results can be be obtained when using more advanced techniques 
such as $l_1$-wavelet regularization, i.e. $R(x) = \|\operatorname{DWT} x\|_1$,  
total variation, or other edge-preserving penalties~\cite{Raj07,Block07,Uecker08b,Liu09}.

\subsection{Discretization}

Numerical reconstruction methods make use of 
discretization, i.e. the unknown image is expanded into 
a sum of basis functions:
\begin{align}
	m(\vec r) \approx \sum_l a_l g_l(\vec r)
\end{align}
Although the choice of this basis has subtle implications 
for results and interpretation, this topic has not drawn
much attention.\footnote{Discussions in terms of ``ideal voxel
functions'' (or ``target voxel shapes'') can be found in earlier 
works~\cite{Pruessmann99,Sodickson01}.}
Most image-domain formulations
based on SENSE use a grid of Dirac pulses to represent the image,
because multiplication with the coil sensitivities in the forward
model is then simply a point-wise multiplication. In contrast, 
k-space methods such as SMASH use a finite Fourier basis.

Discretization has a regularizing effect, i.e. the discretized
problem might have better condition than the continuous problem.
In parallel imaging, this effect can often be seen in the area
outside of the sampled k-space region. Extrapolation to these 
area causes high noise amplification~\cite{Athalye13}. A
discretization scheme which excludes these degrees of freedom will 
be less affected by noise. On the other hand, a small basis
leads to discretization errors because the solution can not
be represented accurately. Both problems can be avoided by using 
fine discretization with a large number of basis function and
explicit regularization to control noise amplification~\cite{Tsao03}.

It should be noted that for parallel imaging the ideal continuous 
solution of Eq.~\ref{eq:var} can usually be computed almost 
perfectly~\cite{Uecker09b}.
Coil sensitivities are very smooth and can be approximated with a 
small number of Fourier coefficients (on an over-sampled FOV). 
The forward operator can then be understood as a convolution
of the Fourier series of the image with a short filter.
Because the acquired k-space data consists of a finite number
of samples, also only a finite number of low-order
Fourier coefficients from the infinite number of coefficients in the 
Fourier series of the image actually appear in the result of this
convolution. For quadratic regularization, a minimum-norm solution is
obtained when the infinite number of remaining higher-order coefficients 
are set to zero. An implementation of the
forward operator requires an a-periodic convolution which can 
be implemented efficiently using a fast Fourier transform (FFT)
algorithm. In practice this differs from a conventional SENSE
implementation only by using zero padding and in the exact 
interpretation of the recovered coefficients.

For non-quadratic regularization, discretization errors
may also arise in the implementation the regularization terms.
In general, oversampling can be used to reduce these errors.
The combination of non-quadratic regularization and oversampling
can also avoid artifacts caused by truncation of the signal 
in the Fourier domain (Gibbs ringing)~\cite{Block08}.
Finally, an important aspect related to discretization is a
common error called an ``inverse crime''~\cite{Colton92}:
When testing a reconstruction algorithms with simulated data,
computing this data using the same discretization scheme as 
used for the reconstruction can result in highly misleading
results. One possibility to avoid this error is the use of 
analytic phantoms~\cite{GuerquinKern12}.

\subsection{Numerical Optimization}

For regular sampling schemes and quadratic regularization,
a solution can be computed directly with matrix inversion,
because the equations decouple into small systems~\cite{Pruessmann99}. 
Although very efficient, this approach is not very flexible. 
Matrix-free iterative methods can be used instead to efficiently 
compute the solution for arbitrary sampling schemes~\cite{Pruessmann01}.
Matrix-free methods are build from procedural implementations 
of the matrix-vector products $y \mapsto F^H y$ and $x \mapsto F x$ 
(or $x \mapsto F^H F x$). For Cartesian sampling, these operations
can be implemented using point-wise multiplications and FFT 
algorithms. For non-Cartesian sampling, efficient non-uniform fast 
Fourier transform (nuFFT) algorithms have been developed to
estimate samples at arbitrary k-space locations~\cite{Sullivan85,Jackson91,Beatty05}.
Even more efficient algorithms can be designed when considering
the combined operator $F^H F$ which appears in the gradient 
of the least-squares data fidelity. For example, the effect of 
sampling in the Fourier domain can be computed exactly 
as a convolution with a truncated point-spread function with
the use of two zero-padded FFTs~\cite{Wajer01}. Overlap-add and 
overlap-save convolution algorithms can be used to exploit the 
compact representation of the coil sensitivities in the 
Fourier domain~\cite{Uecker13b}.

For quadratic regularization, an efficient iterative
algorithm is the conjugate gradient method applied to the normal 
equations~\cite{Pruessmann01}:
\begin{align}
	\left( F^H F + \alpha W \right) x = F^H y~.
\end{align}
It should be noted that the use of a density compensation as 
known from the direct gridding algorithm is neither required 
nor recommended.\footnote{In (non-iterative) gridding, 
the density compensation is a diagonal matrix which approximates 
the inverse of $F^H F$. Combined with the adjoint $F^H$ it 
yields an approximation of the 
pseudo-inverse, i.e. $F^H D \approx F^{\dagger}$.
Including a density compensation into an iterative
optimization method produces solutions different from the 
optimal least-squares solution~\cite{Pruessmann01}.
I.e., naively using $F^H D$ instead of $F^H$ as is sometimes
suggested to improve the condition yields a different
optimization problem:
$ \operatorname{argmin}_x \| \sqrt{D} (F x - y) \|_2^2 + \alpha R(x)$
}

For $l_1$-regularization, the simplest (and slowest)
reconstruction algorithm is iterative soft-thresholding~\cite{Daubechies04}:
\begin{align}
	z_n & = x_n + \alpha F^H \left( y - F x_n \right) \\
	x_n & = T^{-1} \eta_{\lambda}(T z_n)
\end{align}
The first equation is a gradient descent step and
the second update uses soft-thresholding $\eta$ in a
transform basis $T$, e.g. a discrete wavelet transform. This scheme 
converges slowly, but can  be accelerated 
with the addition of a ravine step as in FISTA~\cite{Nesterov83,Beck09}. 
Especially when using multiple convex
penalties $R_n$, a very flexible approach is an extension of
the Alternating Direction Method of Multipliers (ADMM)~\cite{Boyd11,Alfonso11,Ramani11}
that can solve optimization problems of the form
\begin{align}
	\operatorname{argmin}_{x} \sum_{n=1}^L R_n(B_n x)~.
\end{align}
This approach is very flexible and has many advantages
from an implementation point of view, because it splits the
optimization into independent sub-problems. Many different
kinds of regularization terms can easily be integrated
if respective proximal operators of the form
\begin{align}
	prox^R_{\rho}(y) := \operatorname{argmin}_x \frac{\rho}{2} \|x - y\|_2^2 + R(x)
\end{align}
are available in a computationally efficient form. 
For example, the proximal operator for the data fidelity term is 
the $l_2$-regularized least-square inverse which can be computed
efficiently with  the methods of conjugate gradients. The
proximal operator for $l_1$-regularization can be evaluated
simply using soft-thesholding.
Efficient implementations of many advanced algorithms which
make use of parallel programming
can be found Berkeley Advanced Reconstruction Toolbox (BART)~\cite{Uecker14}.

\section{Auto-calibration}

To obtain optimal results in parallel MRI, accurate and up-to-date 
information about the sensitivities of all receive coils is required. 
While approximate coil sensitivities can be computed from the geometry 
of the receive coils using the Bios Savart law, exact sensitivities 
depend  on the loading of the coils and need to be determined with 
high accuracy during the actual measurement.
A pre-scan can provide accurate calibration information,
but this requires that  experimental conditions stay exactly the same for 
the duration of the whole examination. Because this is not always guaranteed,
auto-calibration methods have been developed which perform calibration using a small amount
of additional data acquired during each individual scan~\cite{Jakob98,Griswold02,McKenzie02,Griswold02b}.
Because this reduces overall acceleration, optimal calibration from a minimum 
amount of data is desired. Two advanced techniques are described in the
following: Joint estimation techniques simultaneously estimate image 
content and coil sensitivities from all data, which minimizes
the amount of additional calibration data required. Subspace methods
do not directly estimate sensitivities, but learn a signal subspace 
from calibration data. These algorithms can adapt to experimental conditions 
that violate the sensitivity-based signal model formulated in Equation~\ref{eq:SIG}. 
For this reason, they are more robust to certain kinds of errors.

\subsection{Non-linear Inverse Reconstruction}

Starting with the signal equation (Eq.~\ref{eq:SIG}), but now considering both image
and coil sensitivities as unknowns, one obtains
a non-linear inverse problem related to 
blind multi-channel deconvolution. 
Modelling the
coil sensitivities as smooth functions in a Sobolev space~$H^l(\Omega, \mathbb{C}^N)$,
the non-linear version of the forward operator can be written as:
\begin{align}
	F: L_2(\Omega, \mathbb{C}) \times H^l(\Omega, \mathbb{C}^N) & \rightarrow  C^{\infty}(\mathbb{R}^d, \mathbb{C}^N) \nonumber \\
		x := \left( m, c_1, \cdots, c_N \right) & \mapsto {y}
\end{align}
Many auto-calibrating parallel imaging methods reduce the reconstruction problem 
to a linear problem by first estimating the sensitivities $c_j$ from a subset of 
the data and then solving for the image using these fixed estimates using a
conventional linear reconstruction. Because this is sub-optimal, improved 
algorithms have been developed, which solve the non-linear inverse 
problem~\cite{Bauer05,Ying07,Uecker08}.

In non-linear inversion~\cite{Uecker08}, a regularized solution
is defined as the solution of the minimization problem
\begin{align}
	(\hat m, \hat c_1, \dots, \hat c_N) = \operatorname{argmin}_x \|F x - y\|_2^2 + \alpha R(m) + \beta \sum_{j = 1}^N Q(c_j) ~.
\end{align}
Here, a smoothness penalty $Q(c_j)$ restricts the solution to smooth coil 
sensitivities.
The Iterative Regularized Gauss-Newton Method (IRGNM)~\cite{Bakushinsky93}
is used to iteratively update an estimate of the solution based on
a linearization of the original problem:
\begin{align}
	x_{n+1} - x_n = \operatorname{argmin}_{\Delta x} & \| DF_{x_n} \Delta x +  F x_n - y \|_2^2 \nonumber\\
		&	+ \alpha_n R( \Delta m + m_n ) + \beta \sum_{j = 1}^N Q(\Delta c_j + c_j)
\end{align}
Here, $DF_{x_n}$ is the Frech\'et derivative of $F$ at the current estimate $x_n$
and the regularization parameters $\alpha_n, \beta_n$ are reduced in each iteration step.
The smoothness penalty can be chosen as $Q(c_j) = \|(1 + s\Delta)^l c_j\|_2^2$ (with
some constants $s, l$). This penalty can be transformed into a $l_2$-norm by expressing
the sensitivities using Fourier coefficients re-scaled with a positive definite
diagonal matrix, which avoids bad conditioning of the of the reconstruction problem.
For quadratic regularization of the image, i.e. $R(x) = \|x - x_0\|_2^2$, the algorithm
then has the explicit update rule
\begin{align}
        x_{n+1} - x_n = (DF_{x_n}^{H} DF_{x_n} + \alpha_n I)^{-1} \left( DF_{x_n}^{H} (y - F x_n) +
        \alpha_n( x_n - x_0 ) \right)~.
\end{align}
Non-linear reconstruction methods can be applied to non-Cartesian sampling~\cite{Knoll09,Sheng09,Uecker10}
and extended to include non-linear penalties~\cite{Uecker08b,Knoll12}.

One limitation of non-linear methods is that they may need an initial
guess close to solution to converge to the correct global minimum. 
While it is usually sufficient to set the image to a constant value 
and the coil sensitivities to zero, in some cases a guess closer to 
the true solution is required. In this case, any direct estimation method
can be used to estimate a set of approximate coil sensitivities, 
which can then be used to initialize the non-linear method.

\subsection{Calibration Matrix}

\begin{figure}
\begin{center}
\includegraphics[width=0.9\textwidth]{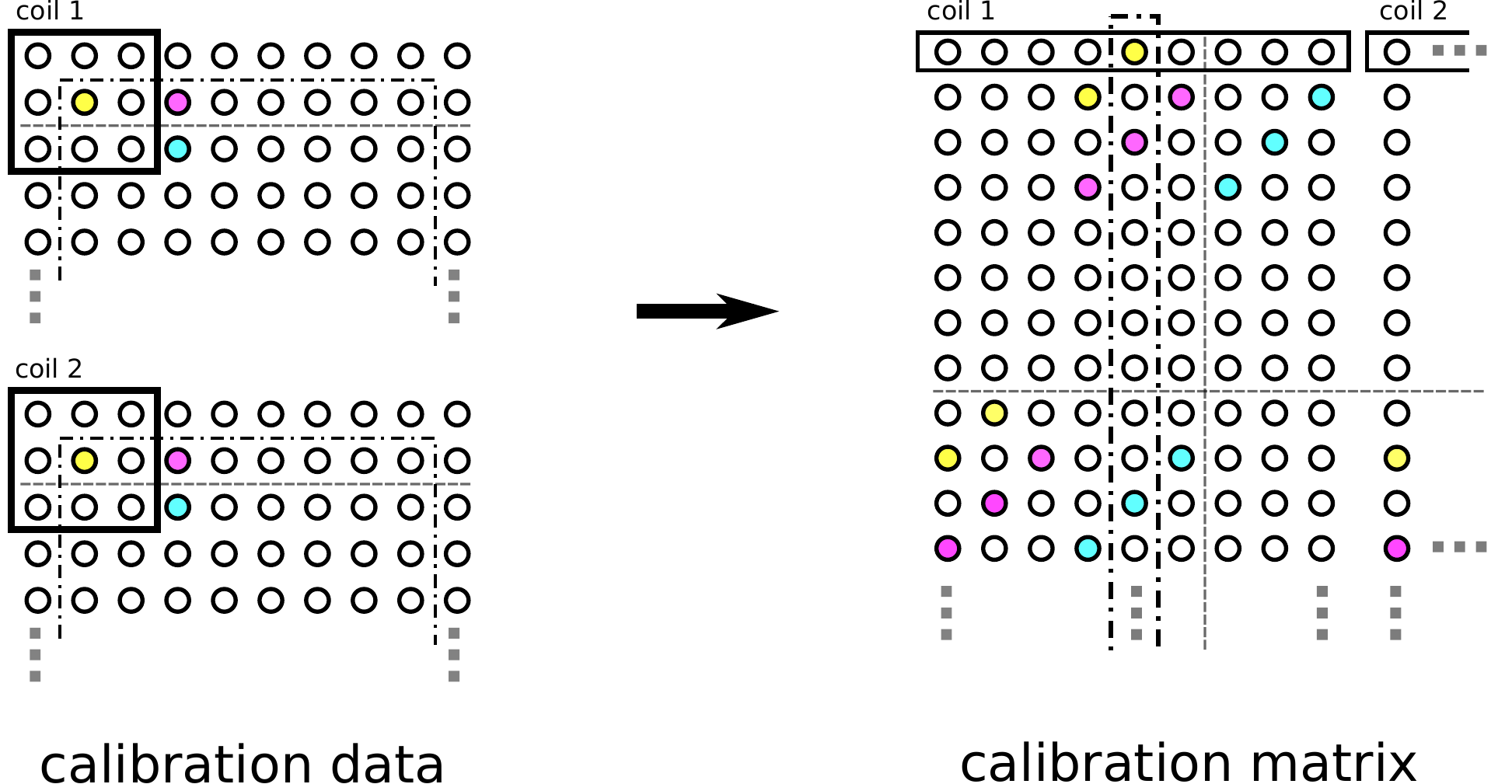}
\end{center}
\caption{Construction of the calibration matrix: Overlapping blocks
of the multi-channel k-space become rows of the calibration 
matrix.}\label{fig:calib}
\end{figure}

\begin{figure}
\begin{center}
\includegraphics[width=0.9\textwidth]{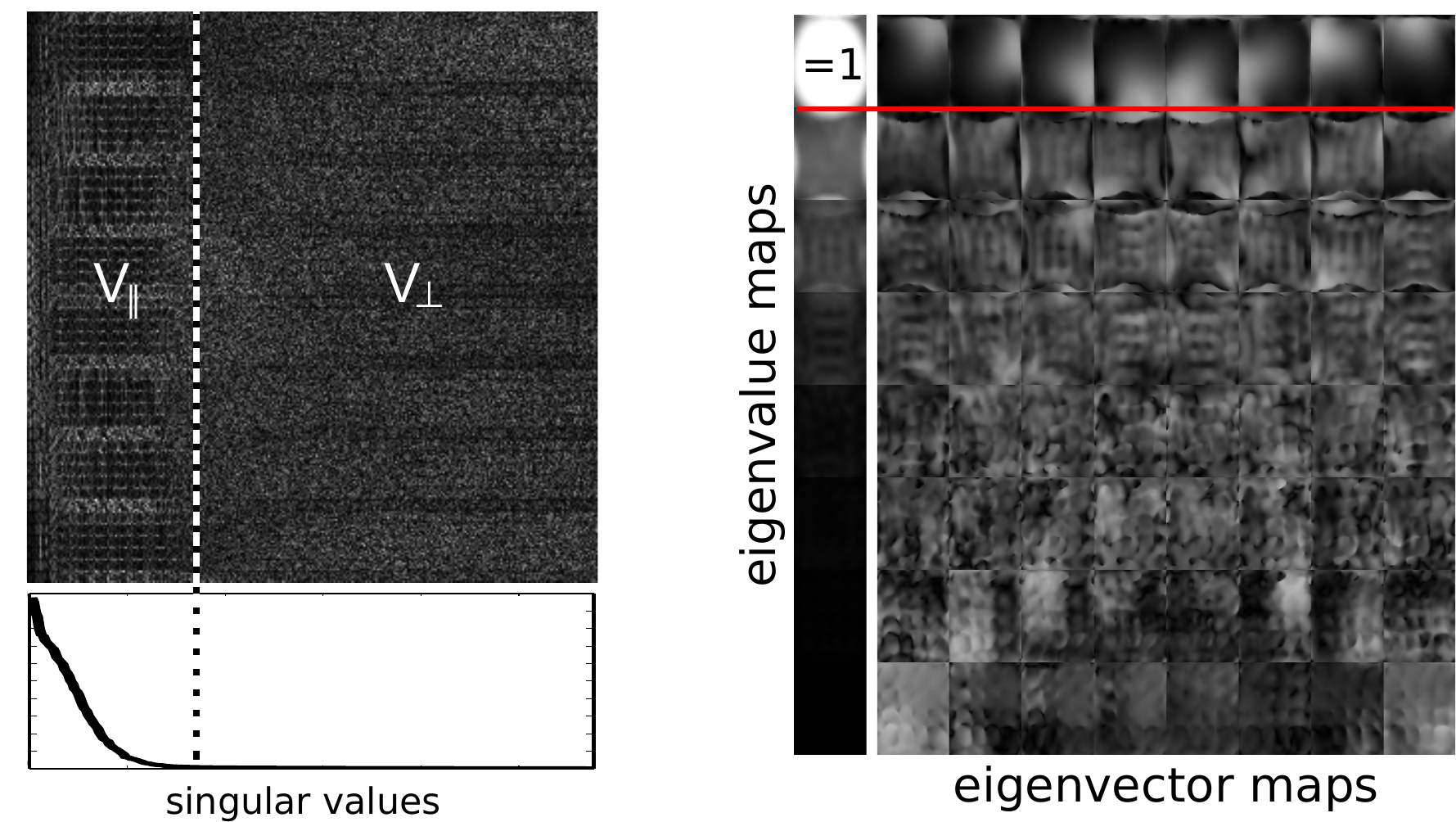}
\end{center}
\caption{Left: Singular-Value-Decomposition (SVD) of the calibration matrix
reveals signal $V_{\parallel}$ and nullspace $V_{\perp}$.
Right: A point-wise eigendecomposition of the operator
$\left[ {\cal F}^{-1} {\cal W} {\cal F} \right]$ is shown.
This operator is derived from the nullspace condition $V_{\perp} R_{\vec k} f = 0$ 
for each overlapping patch in k-space. The sensitivities
(here: from an eight-channel coil) appear at each point
as a eigenvector to eigenvalue one.
}\label{fig:espirit}
\end{figure}

The calibration matrix is a fundamental tool which can
be used to formulate many auto-calibration methods. 
Reconstruction kernels in GRAPPA~\cite{Griswold02} and SPIRiT~\cite{Lustig10}
are null-space vectors of this matrix~\cite{Uecker13}.
The calibration matrix is a multi-dimensional
multi-channel Casorati matrix constructed from fully-sampled
patches in a calibration area in the center of k-space
(see Fig.~\ref{fig:calib}).
It is related to the trajectory matrix from
singular spectrum analysis (SSA)~\cite{Golyandina13},
and also to the lag cross-covariance matrix, which can be
estimated as $\hat \Sigma = \frac{1}{M} C^H C$ 
(with $M$ a normalization constant).
Because coil sensitivities are very smooth, multi-channel
signals have correlations in small local k-space patches. 
This implies that the calibration matrix (and the 
lag cross-covariance matrix) are low-rank, i.e.
have a small signal space and large null space (Fig.~\ref{fig:espirit}). 
If the calibration matrix is constructed from an incomplete
k-space with missing samples, structured low-rank matrix completion 
can be used to recover a completed matrix, which is the
basis of a calibration-less parallel imaging technique known 
as SAKE~\cite{Shin14}.

\subsection{ESPIRiT}

Coil-by-coil reconstruction was originally proposed because
combination of all channels in SMASH-based parallel imaging
sometimes caused phase cancellation~\cite{McKenzie01,Griswold02}. 
In combination with auto-calibration coil-by-coil reconstruction 
has a very advantageous side effect: The reconstruction
becomes robust against certain kinds of inconsistencies - in particular
reconstruction in a tight FOV~\cite{Griswold04}. The fundamental
reason is that the coil-by-coil reconstruction operator does
not enforce the strict signal model of sensitivity-based
reconstruction schemes formulated in Eq.~\ref{eq:SIG}, 
but represents a convex relaxation of this model.

ESPIRiT is a new reconstruction algorithm which exploits this. 
Because of shift-invariance the null-space condition should be
true for all patches in an ideal multi-channel k-space $f$. Let
$R_{\vec k}$ be an operator which extracts a patch around
a given k-space position $\vec k$,
a least-squares version of this condition is then given by
\begin{align}
	\sum_{\vec k} R_{\vec k}^H V_\bot V_\bot^H R_{\vec k}  \, f  & ~=~ 0~.
\end{align}
This can be further transformed to a convolution-type coil-by-coil
operator~$\cal W$, which reproduces ideal k-space signals:
\begin{align}
	\sum_{\vec k} R_{\vec k}^H \left( I - V_\parallel V_\parallel^H \right) R_{\vec k}  \, f & ~=~ 0 \\
	\underbrace{M^{-1} \sum_{\vec k} R_{\vec k}^H V_\parallel V_\parallel^H R_{\vec k}}_{\cal W}  \, f  & ~=~ f\\
	{\cal W} f & ~=~ f
\end{align}
Here, $M$ is the size of a single patch.
Transforming $\cal W$ into the image domain yields
an operator $\left[ {\cal F}^{-1} {\cal W} {\cal F} \right]$
which operates point-wise. Because ${\cal W}$ reproduces ideal
k-space signals, the image-domain version reproduces
the vector of coil images $c_j m$ at each point $\vec r$, i.e.
$\left[ {\cal F}^{-1} {\cal W} {\cal F} \right] c_j(\vec r) m(\vec r) = c_j(\vec r) m(\vec r)$.
In other words, everywhere where the image is non-zero the vector of
sensitivities is a point-wise eigenvector to the 
eigenvalue $1$ of the operator $\left[ {\cal F}^{-1} {\cal W} {\cal F} \right]$.
The eigenvector and eigenvalue maps from a point-wise eigendecomposition
are shown in Figure~\ref{fig:espirit}.
Together, these steps form a computational method to extract accurate coil
sensitivities from the nullspace of the calibration matrix.

If the k-space is corrupted and does not fit the ideal
model, multiple sets of sensitivities can appear in affected image
regions as multiple eigenvectors to eigenvalue one. An extended
forward model can take this additional information
into account. Respective methods offer robustness
to certain kinds of errors similar to auto-calibrating coil-by-coil 
methods such as GRAPPA~\cite{Uecker13}.

\section{Sampling and Reconstruction in k-Space}

\begin{figure}
\begin{center}
\includegraphics[width=0.9\textwidth]{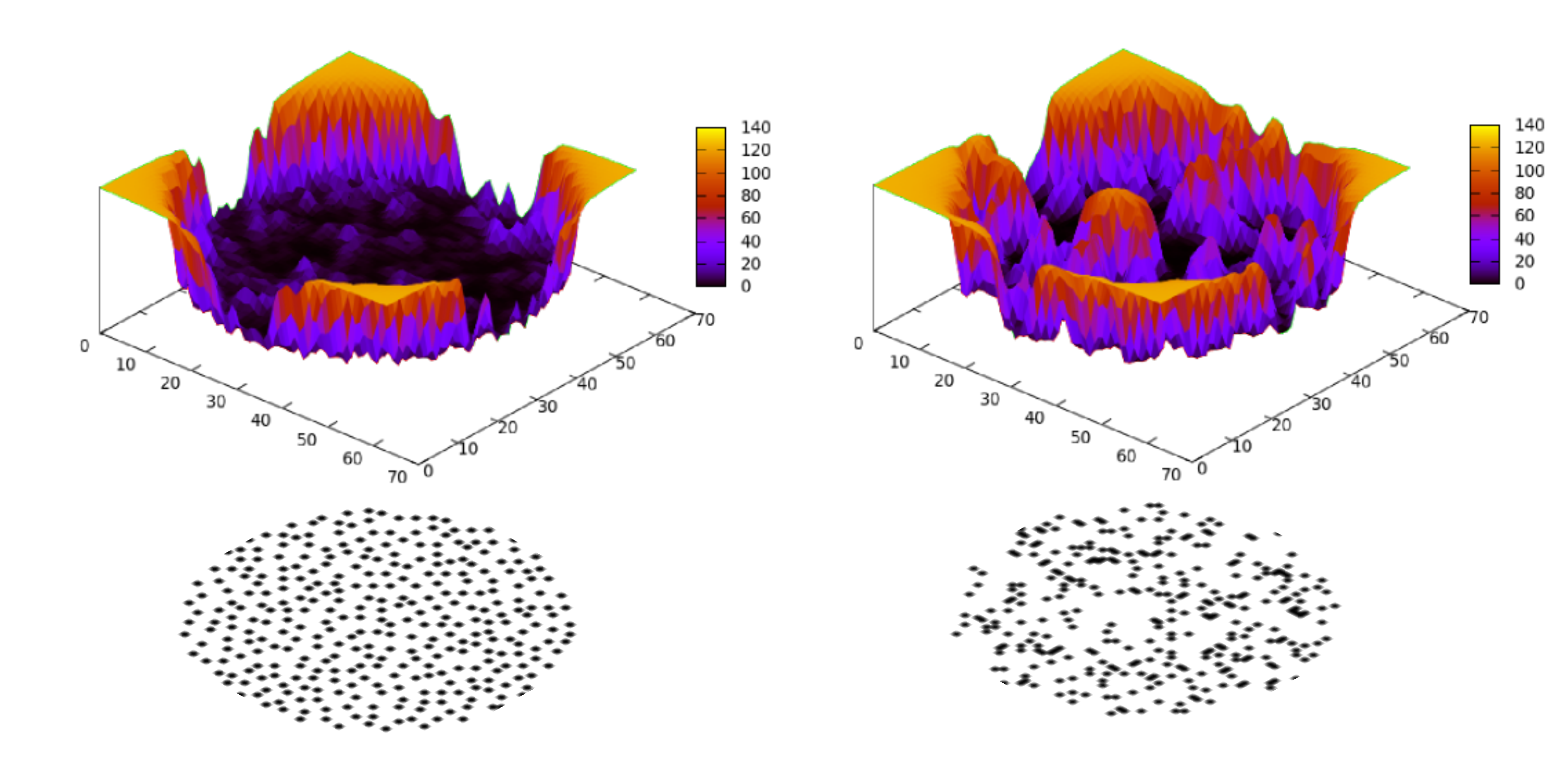}
\end{center}
\caption{The power function for Poisson-disc and random sampling 
for a particular set of coil sensitivities. Outside of the sampled 
area and where random sampling leaves holes without samples the 
power function is high indicating that recovery of the
k-space at this location has a high error when using 
parallel imaging.}\label{fig:POWER}
\end{figure}

While the formulation of parallel imaging as an inverse problem
is a powerful conceptual framework, additional theoretical tools
are required to understand and evaluate different sampling 
schemes in k-space. For this purpose, a formulation 
of parallel imaging as approximation in a Reproducing 
Kernel Hilbert Space (RKHS) has recently been developed~\cite{Athalye13}.
A RKHS is a Hilbert space of functions with continuous (bounded)
point-evaluation functionals. This condition guarantees that sampling
is compatible with the norm of the Hilbert space. 
This formulation of parallel imaging 
yields a unified framework to formulate and
analyze sampling and reconstruction in k-space.

By identifying the norm of the signals in k-space
with the norm of the corresponding images in $L(\Omega, \mathbb{C})$,
it can be shown that the ideal k-space signals are 
an RKHS with the matrix-valued kernel
\begin{align}
	K_{st}(\vec k, \vec l) 
				& = \int_{\Omega} \textrm{d}\vec r \, c_s(\vec r) \overline{c_t(\vec r)} 
					e^{-2\pi i \vec r \cdot (\vec k - \vec l)}~. \label{eq:kern}
\end{align}
This kernel captures all local correlations in k-space induced by 
the sensitivities and exploited in parallel imaging algorithms
for recovery of missing k-space samples. Given this kernel, a standard 
formula from  approximation theory
can be applied to obtain interpolation coefficients $u_j^{s, \vec k}$
for all channels $s = 1, \dots, N$ and known k-space samples $\vec k \in S$ 
for interpolation to arbitrary k-space positions:
\begin{align}
	\sum_{s = 1}^N \sum_{\vec k \in S} K_{st}(\vec k, \vec l) u_j^{s, \vec k}(\cdot) = K_{jt}(\cdot, \vec l)
\end{align}
With these interpolation coefficients, unknown values in k-space can then 
be recovered from the acquired samples $f_j(\vec k(t_l))$ with the 
interpolation formula
\begin{align}
	{\hat f}_j(\vec k) = \sum_{t = 1}^N \sum_{\vec l \in S} f_t(\vec l) {u_j}^{t, \vec l}(\vec k)~.
\end{align}
When no regularization is used in the computation of coefficients
the recovered ideal k-space corresponds to the best approximate solution 
defined before. The interpolation formulas used in GRAPPA 
and SPIRiT and similar methods are local variants of this
formula with an empirical estimate of the
ideal kernel (Eq.~\ref{eq:kern}).
In addition to this interpolation formula, the link to approximation 
theory yields new insights into sampling in k-space.
In particular, a point-wise error bound in k-space can be derived~\cite{Schaback93b}:
\begin{align}
	 | {f_j}(\vec k) - {\hat f_j}(\vec k) |^2 \leq \| {f} \|^2 \cdot \underbrace{\left( K_{jj}(\vec k, \vec k) - \sum_{t = 1}^N \sum_{\vec l \in S} K_{jt}(\vec k, \vec l) \overline{u_j^{t, \vec l}(\vec k)} \right)}_{P^2_n(\vec k)}
\end{align}
The power function $P_j$ is computed from the kernel $K$
and the interpolation functions $u$ and depends only on the coil
sensitivities and the sample locations. It can be used to analyze 
the properties of different sampling schemes for parallel imaging 
independent from any actual imaging data. 
Figure~\ref{fig:POWER} shows the Power function 
for two different sampling schemes computed for 
a particular set of coil sensitivities.


\section{Compressed Sensing Parallel Imaging}


Compressed sensing is based on the idea that randomized 
under-sampling schemes produce incoherent noise-like artifacts
in a transform domain which can then be suppressed 
using denoising to iteratively recover
the original signal~\cite{Donoho06,Candes06}.
It exploits the compressibility of the image 
information, i.e. a sparse representation in a transform
domain, to make the sparse signal coefficients stand out 
from the incoherent noise.
Non-linear regularization terms can then be used to 
suppress the incoherent artifacts and recover a sparse
representation of the image from the under-sampled data.
Because MRI acquires data in the Fourier domain and 
is flexible enough to use almost arbitrary sampling schemes,
this idea can be applied directly~\cite{Block07,Lustig07}.

\begin{figure}
\begin{center}
\includegraphics[width=0.5\textwidth]{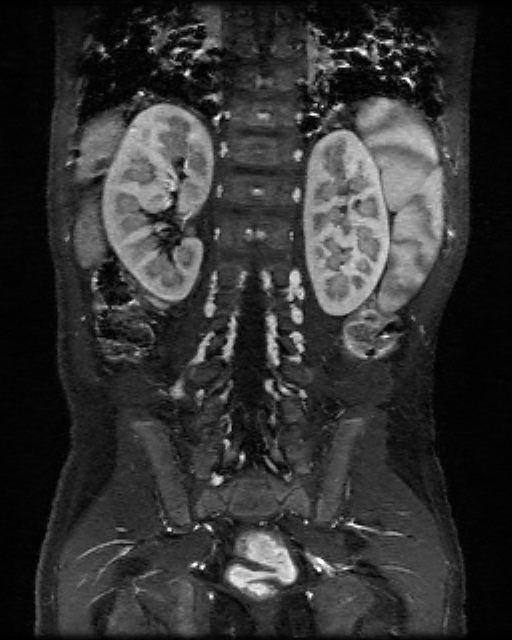}
\end{center}
\caption{Parallel imaging compressed sensing
reconstruction of an abdomen from a pediatric patient
using $l_1$-ESPIRiT. The acquisition made us of 
variable-density Poisson-disc sampling with an acceleration 
factor of seven.}\label{fig:CSPI}
\end{figure}

Parallel imaging can be synergistic-ally combined with
compressed sensing~\cite{Block07,Liu09,Lustig10,Otazo10}.
This combination leads to exactly the same optimization 
problems already considered for parallel imaging alone, 
but requires incoherent sampling schemes 
suitable for compressed sensing. The most important 
schemes in practical use are variable-density Poisson-disc 
sampling and radial trajectories.
Poisson-disc sampling guarantees that samples are not too
close  together. This would waste sampling time, because 
k-space positions which are close are highly correlated
and can already be recovered using parallel imaging. 
Variable-density schemes have several advantages: They
equalize the power  spectrum of the missing samples, 
provide graceful degradation in case full recovery is 
not possible, and can be used for auto-calibrating 
parallel imaging when the k-space center is 
fully sampled.
Methods which combine parallel imaging and compressed
sensing represent the state of the art in image
reconstruction, as demonstrated by their use 
in demanding applications such as in 
pediatric imaging without
sedation~\cite{Vasanawala10,Zhang14,Cheng14}.
Figure~\ref{fig:CSPI} shows an
image from a pediatric patient reconstructed
with parallel imaging compressed sensing 
at an acceleration factor of about seven.

\section{Conclusion}

Image reconstruction for parallel imaging can be formulated
as an inverse problem. Based on this formulation, advanced 
iterative algorithms can be developed which (1) make use
of optimal (Cartesian or non-Cartesian) sampling schemes, and
(2) extend parallel imaging with advanced non-linear regularization
terms. These ideas are combined in recent methods for 
compressed sensing parallel imaging, which currently 
represent the state of the art in image reconstruction.

\bibliography{pics}

\end{document}